# Proof of the Irrationality of the Square Root of Two in Babylonian Geometry Problem Tablets

**Abstract**
One of the greatest achievements of Greek mathematics is the proof that the square root of 2 is irrational. It has not been thought that the Babylonians appreciated the concept of irrationality and certainly that they did not prove that the square root of two is irrational. Here we show that two Babylonian geometry problem tablets contain a simple proof of the irrationality of the square root of two. It is not known, as yet, if the Babylonians appreciated that these tablets indeed contained this proof.

One of the greatest achievements of Greek mathematics is the proof that the square root of 2 is irrational [1]. The proof is traditionally credited to the circle of Pythagoras (c. 570-c. 495 BCE) [2], however this specific attribution is disputed [3]. This Greek proof is algebraic and proceeds by contradiction: assume that the square root of two is rational so therefore it can be written as the quotient of two integers (p/q) with this fraction in lowest terms. Simple algebraic manipulations then quickly yield a contradiction as both p and q are found to be even so the quotient is not in lowest terms. A century after Pythagoras Theodorus of Cyrene (5$^{th}$ century BCE) was able to prove the irrationality of the square root of three, five and other numbers up to seventeen [4]. Despite the Babylonians' prowess in mathematics and astronomy and indeed an Old Babylonian tablet (known as YBC 7289 [5-8] (c. 1800-1600 BCE) containing a diagram of a square with diagonals and inscribed in cuneiform the best three digit sexagesimal approximations (accurate to six decimal digits) of both the square root of two and its reciprocal, it has not been thought by mathematicians [5, 7, 8] or historians of mathematics [1, 6, 9, 10] that the Babylonians knew, or certainly did not prove that the square root of two is irrational. Here we point out that a problem on the Old Babylonian tablet BM 15285 [11] contains a simple geometric proof that the square root of two is irrational.

Problem twelve (of a presumed 41 geometric area problems) on BM 15285 states: The side of the square is 60 rods. Inside it [I drew]16 wedges [triangles]. What are their areas? The figure for the problem is given below (Fig. 1). First we show that the tablet contains a geometric proof of what we call today the Pythagorean theorem for the special case of an isosceles right triangle. The problem figure can be constructed with a compass and straight edge: Draw the large square by constructing lines perpendicular to one edge and marking off on these lines segments of length of the first edge. Let us call the four vertices of the large square A, B, C and D. Draw the two diagonals and call the midpoint (by symmetry) O. On segments AO, BO, CO and DO find midpoints L, M, N, and R. Draw lines through L, M, N, and R perpendicular to BD, AC, BD and AC respectively and call the intersection of these lines with the sides of the large square W, X, Y and Z. Finally connect points L, M, N, and R. By symmetry all sixteen triangles have equal areas that we denote as T. Consider now, for example, the isosceles right triangle (by the above construction) LMN. By considering square LMNR $(LM)^2 = (MN)^2 = 4T$. By considering square WXYZ it is clear that $(LN)^2 = 8T$. Thus $(LM)^2 + (MN)^2 = (LN)^2$ and the Pythagorean theorem for isosceles triangles is established within the problem figure.

Another way to write this theorem is that if we denote by H the length of the hypotenuse of an isosceles right triangle and S by the length of a side $H^2 = 2S^2$. From this we also see immediately that the hypotenuse of an isosceles right triangle must be longer than a side.

Now to establish that the problem figure also contains a proof of the irrationality of the square root of two. Consider the large square and assume to the contrary that √2 is rational. Then denoting by S the length of the side of the large square and H the length of the main diagonal we have since $H^2=2S^2$ that H/S (=√2) is a rational number (a quotient of two integers) which has been reduced to lowest terms. But now consider triangle AOB, itself by construction also an isosceles right triangle. The hypotenuse of triangle AOB is AB that has length S. The sides of AOB (AO and BO) are by symmetry length H/2. Since $H^2$ is equal to twice an integer ($S^2$) $H^2$ is even. But if the integer $H^2$ is even, H also must be even. Thus H/2 is an integer. By the Pythagorean theorem proved above S/(H/2)= √2. But since S and H/2 are integers and S<H and H/2<S there is a contradiction that H/S was reduced to lowest terms. Thus √2 is irrational. QED.

If one assumes the Pythagorean theorem the construction on YBC 7289 also contains the same irrationality proof. A similarity of the figure on YBC 7289 to some of the figures on BM 15285 has been noted [7].

It is not yet known if the Babylonians appreciated that these problem tablets contain a proof of the irrationality of the square root of two, or even if they understood the idea of an irrational number [1, 2, 5-11]. Perhaps presciently Neugebauer noted [1], "But all the foundations were laid which could have given this result [irrationality of √2] to a Babylonian mathematician …." As we now see if fact the Babylonian geometric problem tablets indeed contain a proof of this fact! Further investigation is merited to see if indeed the Babylonians appreciated and understood the implications of the proof contained in their tablets and thus that the Babylonian geometers predated the Greeks in proving this profound mathematical fact. Regardless, the problem xii is a "Book Proof" indeed, a beautiful, simple and without words "Tablet Proof."

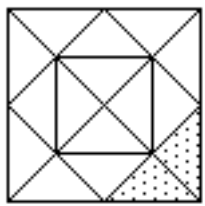

**Figure 1.** Figure for problem xii on tablet BM 15285 [11]. The speckled area indicates a broken off part of the tablet.

Benjamin M. Altschuler[1], Eric L. Altschuler, MD, PhD[2]*

[1]The Fieldston School, Bronx, NY, USA

[2]Lewis Katz School of Medicine at Temple University, Philadelphia, PA, 19140, USA, eric.altschuler@temple.edu